\documentclass[11pt]{amsart}
\usepackage{amsmath}
\usepackage{amssymb}
\usepackage{amsthm}
\usepackage{latexsym}
\date{}

\def\opn#1#2{\def#1{\operatorname{#2}}} 
\opn\chara{char} \opn\length{\ell} \opn\pd{pd} \opn\rk{rk}
\opn\projdim{proj\,dim} \opn\injdim{inj\,dim} \opn\rank{rank}
\opn\depth{depth} \opn\codepth{codepth} \opn\grade{grade}
\opn\height{height} \opn\embdim{emb\,dim} \opn\codim{codim}

\opn\Tr{Tr} \opn\bigrank{big\,rank}
\opn\superheight{superheight}\opn\lcm{lcm}
\opn\trdeg{tr\,deg}%
\opn\reg{reg} \opn\lreg{lreg} \opn\skel{skel} \opn\Gr{Gr}
\opn\dim{dim} \opn\arithdeg{arithdeg}

\opn\div{div} \opn\Div{Div} \opn\cl{cl} \opn\Cl{Cl}
\opn\Spec{Spec} \opn\Supp{Supp} \opn\supp{supp} \opn\Sing{Sing}
\opn\Ass{Ass}
\opn\Ann{Ann} \opn\Rad{Rad} \opn\Soc{Soc}
\opn\Sym{Sym} \opn\Ker{Ker} \opn\Coker{Coker} \opn\Im{Im}
\opn\Hom{Hom} \opn\Tor{Tor} \opn\Ext{Ext} \opn\End{End}
\opn\Aut{Aut} \opn\id{id} \opn\ini{in} \opn\tr{tr}

\opn\nat{nat}\opn\it{it}
\opn\pff{proof}
\opn\Pf{proof} \opn\GL{GL} \opn\SL{SL} \opn\mod{mod} \opn\ord{ord}
\opn\aff{aff} \opn\con{conv} \opn\relint{relint} \opn\st{st}
\opn\lk{lk} \opn\cn{cn} \opn\core{core} \opn\vol{vol}
\opn\link{link} \opn\star{star} \opn\skel{skel}
\opn\gr{gr}

\def\pot#1#2{#1[\kern-0.28ex[#2]\kern-0.28ex]}

\opn\dirlim{\underrightarrow{\lim}}
\opn\inivlim{\underleftarrow{\lim}}
\def\Implies{\ifmmode\Longrightarrow \else
     \unskip${}\Longrightarrow{}$\ignorespaces\fi}
\def\implies{\ifmmode\Rightarrow \else
     \unskip${}\Rightarrow{}$\ignorespaces\fi}
\def\iff{\ifmmode\Longleftrightarrow \else
     \unskip${}\Longleftrightarrow{}$\ignorespaces\fi}

\let\:=\colon
\newtheorem{Theorem}{Theorem}[section]
\newtheorem{Lemma}[Theorem]{Lemma}
\newtheorem{Corollary}[Theorem]{Corollary}
\newtheorem{Proposition}[Theorem]{Proposition}
\newtheorem{Remark}[Theorem]{Remark}

\newtheorem{Example}[Theorem]{Example}

\let\epsilon\varepsilon
\let\phi=\varphi
\let\kappa=\varkappa
\def\qed{\ifhmode\textqed\fi
   \ifmmode\ifinner\quad\qedsymbol\else\dispqed\fi\fi}
\def\textqed{\unskip\nobreak\penalty50
    \hskip2em\hbox{}\nobreak\hfil\qedsymbol
    \parfillskip=0pt \finalhyphendemerits=0}
\def\dispqed{\rlap{\qquad\qedsymbol}}

\opn\Gin{Gin}

\opn\inii{in} \opn\inim{inm} \opn\rate{rate}

\numberwithin{equation}{section}

\begin{document}

\title[Cohen-Macaulay edge ideal]{Cohen-Macaulay edge ideal whose height is half of the number of vertices}
\author[Marilena Crupi]{Marilena Crupi}
\address[Marilena Crupi]{Dipartimento di Matematica,
Universita' di Messina, Salita Sperone, 31, 98166 Messina, Italy. Fax number: +39 090 393502}
\email{mcrupi@unime.it}
\author[Giancarlo Rinaldo]{Giancarlo Rinaldo}
\address[Giancarlo Rinaldo]{Dipartimento di Matematica,
Universita' di Messina, Salita Sperone, 31, 98166 Messina, Italy. Fax number: +39 090 393502}
\email{rinaldo@dipmat.unime.it}
\author[Naoki Terai]{Naoki Terai}
\address[Naoki Terai]{Department of Mathematics, Faculty of Culture
and Education, Saga University, Saga 840--8502, Japan. Fax number: }
\email{terai@cc.saga-u.ac.jp}

\date{\today}

\begin{abstract}
We consider  a class of graphs $G$ such that  the height of the edge ideal $I(G)$
is half of the number $\sharp V(G)$ of the vertices.
We give  Cohen-Macaulay criteria for such  graphs.
\end{abstract}

\maketitle

\section*{Introduction}

In this article a graph means a simple graph without loops and
multiple edges. Let $G$ be a graph with the vertex set $V(G) = \{x_1,
\dots , x_n\}$ and with the edge set $E(G)$. Let $S = K[x_1, \ldots,
x_n]$  be the polynomial ring in $n$ variables over a field $K$. The
\textit{edge ideal} $I(G)$, associated to $G$, is the ideal of $S$
generated by the set of all squarefree monomials $x_ix_j$ so that
$x_i$ is adjacent to $x_j$. For this ideal the following theorem
\cite {GV1} is known:
\medskip
\begin{Theorem}
Suppose $G$ is an unmixed graph without isolated vertices. Then we
have $2 \height I(G) \ge \sharp V(G)$.
\end{Theorem}
\bigskip In this paper we treat the class of graphs for which  the
above equality holds, i.e., we consider an unmixed graph without
isolated vertex with $2 \height I(G) = \sharp V(G)$. Such a class
of graphs is rich, because it includes all the unmixed bipartite
graphs and all the grafted graphs. Herzog-Hibi \cite {HH} gave
beautiful theorems on Cohen-Macaulay edge ideals of bipartite
graphs. Our purpose in this article is to generalize their results
for our class of graphs.

It is known that a graph $G$ in our class has a perfect
matching, we may assume
that

\vspace{.2in}
(*)
$V(G)= X \cup Y$, $X \cap Y = \emptyset$, where $X=\{x_1, \ldots, x_n\}$ is a minimal vertex cover of $G$
and $Y=\{ y_1, \dots, y_n\}$ is a maximal independent set of $G$ such that $\{x_1y_1, \dots, x_ny_n\} \subset E(G)$.
\vspace{.2in}

Hence $\{ x_1-y_1, \ldots, x_n- y_n$\} is a system of parameters
of $S/I(G)$. In Sections 3 and 4, using this, we give the
following characterization of Cohen-Macaulayness, which is
similar to the case of bipartite graph (see \cite {HH}).

\begin{Theorem}
Let $G$ be an unmixed graph with $2n$ vertices, which are not isolated,
and with $\height I(G) = n$.
Then the following conditions are equivalent:
\begin{enumerate}
 \item $G$ is Cohen-Macaulay.
 \item $\Delta(G)$ is strongly connected.
 \item There is a unique perfect matching in $G$.
 \item $\Delta(G)$ is shellable.
\end{enumerate}
\end{Theorem}

Note that it includes equivalence between Cohen-Macaulayness and shellability
as in the bipartite graphs (see \cite{EV}).

We also have a Cohen-Macaulay criterion which is similar to
Herzog-Hibi (\cite{HH}, Theorem 3.4):

\begin{Theorem} Let $G$ be a graph with $2n$ vertices, which are not isolated, and with $\height I(G) =n$.
We assume the conditions (*) and
\begin{center}
(**) $x_iy_j \in E(G)$ implies $i \leq j$.
\end{center}
Then the following conditions are equivalent:
\begin{enumerate}
\item $G$ is Cohen-Macaulay.
\item $G$ is unmixed.
\item The following conditions hold:
\begin{enumerate}
\item[(i)] If $z_ix_j, y_jx_k \in E(G)$, then $z_ix_k \in E(G)$ for distinct $i, j, k$ and
for $z_i \in \{ x_i, y_i \}$.
\item[(ii)] If $x_iy_j \in E(G)$, then $x_ix_j \notin E(G)$.
\end{enumerate}
\end{enumerate}
\end{Theorem}

Although in Herzog-Hibi \cite {HH}  Alexander duality plays an
important role for their proof, we give a direct and elementary
proof without it. The Herzog-Hibi criterion for bipartite graphs
was discussed by many authors in literature that gave alternative
proofs for it (see \cite{HY}, \cite{TuVi})

In Section 5 we introduce a new class of graphs which we call
B-grafted graphs. They are a generalization of  grafted graphs
introduced by Faridi \cite{SF}. If $G$ is  an unmixed B-grafted
graph, then we have  $2 \height I(G) = \sharp V(G)$. Hence
applying our main result, we show:

\begin{Theorem} The B-grafted graph $G(H_0; B_1, \ldots, B_p)$ is Cohen-Macaulay (unmixed, respectively) if and
only if every bipartite graph $B_i$ is   Cohen-Macaulay (unmixed, respectively) for $i = 1, \ldots, p$.
\end{Theorem}

See Sections 1 and 5 for undefined concepts and notation.

\section{Preliminaries}
In this section we recall some concepts and a notation on graphs and on simplicial
complexes that we will use in the article.

Let $G$ be a graph with the vertex set $V(G) = \{x_1\dots, x_n\}$
and with the edge set $E(G)$.
The \textit{induced subgragh}  $G\vert _{W}$ by $W\subset V(G)$  is defined by
$$
G\vert _{W}=(W, \{ e \in E(G); e \subset W \}).
$$
For $W \subset V(G)$ we denote $G \vert _{V(G)\setminus W}$ by $G-W$.
For $F \subset E(G)$ we denote $(V(G), E(G) \setminus F)$ by $G-F$.
For a family $F$ of 2-element subsets of $V(G)$ we denote $(V(G), E(G) \cup F)$ by $G+F$.

A subset $C \subset V(G)$ is a \textit{vertex cover} of $G$ if
every edge of $G$ is incident with at least one vertex in $C$.
A  vertex cover $C$ of $G$ is called \textit{minimal} if there
is no proper subset of $C$ which is a vertex cover of $G$.
A subset $A$ of $V(G)$ is called an \textit{independent set} of $G$
if no two vertices of $A$ are adjacent.
An independent set $A$ of $G$ is \textit{ maximal } if there exists no independent set
which properly includes $A$.
Observe that $C$ is a minimal vertex cover of $G$ if and only if $V(G) \setminus C$
is a maximal independent set of $G$.
And also note that $\height I(G)$ is equal to the smallest number $\sharp C$ of vertices among
all the minimal vertex covers $C$ of $G$.
A graph $G$ is called \textit{unmixed} if all the minimal vertex
covers of $G$ have the same number of elements. A graph $G$ is
called \textit{Cohen-Macaulay} if $S/I(G)$ is a Cohen-Macaulay
ring, where $S = K[x_1, \dots, x_n]$ is a polynomial ring for a
field $K$. Refer \cite{Die}, \cite{Vi1} for detailed information
on this subject.

\vspace {.2in}

Set $V = \{x_1, \ldots, x_n\}$. A \textit{simplicial complex}
$\Delta$ on the vertex set $V$ is a collection of subsets of $V$
such that (i) $\{x_i\} \in \Delta$  for all $x_i \in V$ and (ii)
$F \in \Delta$ and $G\subseteq F$ imply $G \in \Delta$. An
element $F \in \Delta$ is called a \textit{face} of $\Delta$. For
$F \subset V$ we define the \textit{dimension} of $F$ by $\dim F
= \sharp F -1$, where $\sharp F$ is the cardinality of the set
$F$. A maximal face of $\Delta$  with respect to inclusion is
called a \textit{facet} of $\Delta$.
If all facets of $\Delta$ have the same
dimension, then $\Delta$ is called \textit{pure}.

A pure simplicial complex $\Delta$ is called \textit{shellable} if
the facets of $\Delta$ can be given a linear order $F_1, \ldots,
F_m$ such that for all $1 \leq j < i \leq m$, there exist some $v
\in F_i \setminus F_j$ and some $k \in \{1, \dots, i-1\}$ with $F_i
\setminus F_k = \{v\}$.

Moreover, a pure simplicial complex $\Delta$ is \textit{strongly
connected} if for every two facets $F$ and $G$ of $\Delta$  there
is a sequence of facets $F = F_0, F_1, \ldots, F_m = G$ such that
$ \dim ( F_i \cap F_{i+1}) = \dim \Delta -1$ for each $i =0, \ldots,
m-1$.

If $G$ is a graph, we define the \textit{complementary simplicial
complex} of $G$ by
\[\Delta(G) = \{A \subseteq V(G) \, : \, \mbox{$A$ is an independent set
in}\, G\} .\] Observe that $\Delta(G)$ is the Stanley-Reisner
simplicial complex of $I(G)$.

\section{Unmixedness}

In this section we survey unmixed graphs whose edge ideals
have the height that is half of the number of vertices.

\begin{Lemma} \label{Pro:main} Let $G$ be an unmixed graph  with non-isolated $2n$
vertices and with $\height I(G)= n$.
Then $G$ has a perfect matching.
\end{Lemma}

The proof is clear from (\cite {GV2}, Remark 2.2). By the lemma  for
an unmixed graph $G$  with $2n$ vertices, which are not isolated,
and with $\height I(G)=n$, we may assume

\medskip

(*)
 $V(G)=X \cup Y$, $X \cap Y=\emptyset$,
 where $X=\{x_1,\ldots,x_n\}$ is a minimal vertex cover of $G$
 and $Y=\{y_1,\ldots, y_n\}$ is a maximal independent set of $G$
 such that  $\{x_1 y_1, \ldots, x_ny_n\} \subset E(G)$.

\medskip
From now on, set $S=K[x_1,\ldots,x_n, y_1,\ldots, y_n]$ for a field $K$
and $I(G)$ is an ideal of $S$.
By Lemma 2.1 we have the following ring-theoretic properties of $S/I(G)$:

\begin{Corollary}\label{pro:covxy}
Let $G$ be an unmixed graph with $2n$ vertices, which are not
isolated, and with $\height I(G)=n$. We assume the condition (*).
Then \begin{enumerate}
\item[(i)] Each minimal prime
ideal of $I(G)$ is of the form
\[
 (x_{i_1},\ldots,x_{i_k},y_{i_{k+1}},\ldots,y_{i_n}),
\]
where $\{i_1,\ldots,i_n\}=\{1,\ldots,n \}$.
\item[(ii)] $\{x_1-y_1, \ldots, x_n- y_n \}$ is a system of parameters of
$S/I(G)$.
\end{enumerate}
\end{Corollary}

%
%

For later use we give a characterization of  the unmixedness for our
graphs, that is a more detailed description, but  a  special case of
a more general result (\cite{MoReVi}, Theorem2.9):

\begin{Proposition}\label{th:unmix_ijk}
Let $G$ be a graph with $2n$ vertices, which are not isolated,  and  with
$\height I(G)=n$. We assume the condition (*). Then $G$ is
unmixed  if and only if the following conditions hold:
 \begin{enumerate}
   \item[(i)] If $z_i x_j$, $y_j x_k\in E(G)$ then $z_i x_k\in E(G)$
    for distinct $i$, $j$ and $k$ and for $z_i\in \{ x_i, y_i\}$.
   \item[(ii)] If $x_iy_j \in E(G)$ then $x_ix_j \notin E(G)$.
 \end{enumerate}
\end{Proposition}

%
%
%
%
%

\par\bigskip

\section{Cohen-Macaulayness}

In this section we give  combinatorial characterizations of Cohen-Macaulay graphs whose edge ideals
have the height that is half of the number of vertices.

First we introduce an operator that allows us to construct a new graph.
Let $G$ be a graph with $2n$ vertices, which are not isolated,
and with $\height I(G) = n$.
We assume the condition (*).

For any $i \in [n]:= \{1, \dots, n\}$ set
\[
E_i:=\{k \in [n] \,:\, x_ky_i \in E(G)\} \setminus \{i\},
\]
 and define the graph $O_i(G)$ by
\[
O_i(G) := G - \{x_ky_i\,:\, k \in E_i\} + \{x_kx_i\,:\,k \in E_i\}.
\]
Then for every subset $T:=\{i_1, \dots, i_{\ell}\}$ of the set
$[n]$, we define
\[O_T(G) = O_{i_1}O_{i_2}\cdots O_{i_{\ell}}(G).\]

Note that $O_T(G)$ is a graph with $2n$ vertices, which are not isolated,
and with $\height I(G) = n$ satisfying the condition (*).

\begin{Example}
Let $T=\{2,3\}$, then
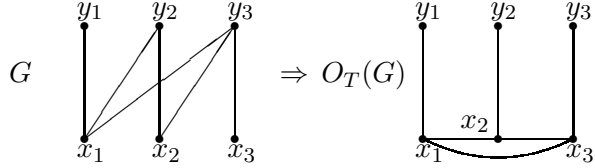
\[
\setlength{\unitlength}{.5cm}
\begin{picture}(15,5)

\put(1,1){\line(0,1){3}}
\put(1,1){\line(2,3){2}}
\put(1,1){\line(4,3){4}}

\put(3,1){\line(0,1){3}}
\put(3,1){\line(2,3){2}}

\put(5,1){\line(0,1){3}}

\put(1,1){\circle*{.2}}
\put(.8,.5){$x_1$}

\put(3,1){\circle*{.2}}
\put(2.8,.5){$x_2$}

\put(5,1){\circle*{.2}}
\put(4.8,.5){$x_3$}

\put(-1,2.5){$G$}

\put(1,4){\circle*{.2}}
\put(.8,4.3){$y_1$}

\put(3,4){\circle*{.2}}
\put(2.8,4.3){$y_2$}

\put(5,4){\circle*{.2}}
\put(4.8,4.3){$y_3$}

\put(6.2,2.5){$ \Rightarrow\, O_T (G)$}

\put(10,1){\line(0,1){3}}
\put(10,1){\line(1,0){2}}

\qbezier(10,1)(12,0)(14,1)

\put(12,1){\line(0,1){3}}
\put(12,1){\line(1,0){2}}

\put(14,1){\line(0,1){3}}

\put(10,1){\circle*{.2}}
\put(9.8,.5){$x_1$}

\put(12,1){\circle*{.2}}
\put(11,1.3){$x_2$}

\put(14,1){\circle*{.2}}
\put(13.8,.5){$x_3$}

\put(10,4){\circle*{.2}}
\put(9.8,4.3){$y_1$}

\put(12,4){\circle*{.2}}
\put(11.8,4.3){$y_2$}

\put(14,4){\circle*{.2}}
\put(13.8,4.3){$y_3$}


\end{picture}
\]
\end{Example}

The next proposition shows that Cohen-Macaulayness of $G$ can be checked by unmixedness of
all the deformations $O_T(G)$ of $G$.

\begin{Proposition}\label{Prop:CM} Let $G$ be an unmixed graph with $2n$ vertices, which are not isolated,
and $\height I(G) = n$. We assume the condition (*). Then the
following conditions are equivalent:
\begin{enumerate}
\item $G$ is Cohen-Macaulay.
\item $O_T(G)$ is Cohen-Macaulay for every subset $T$ of $[n]$.
\item $O_T(G)$ is unmixed for every subset $T$ of $[n]$.
\end{enumerate}
\end{Proposition}

\begin{proof} Set
$S = K[x_1, \dots, x_n, y_1, \dots, y_n]$,
$S_k =K[x_1, \dots, x_n, y_{k+1}, \dots, y_n]$,
and $G_k = O_{T_k}(G)\vert_{X\cup\{y_{k+1}, \dots, y_n\}}$.

(1) $\Longrightarrow$ (2). By relabeling, we may assume that $T=[k]$.
Let $G$ be a Cohen-Macaulay graph.
Then
\[
S/(I(G) + (x_1 -y_1, \ldots , x_k -y_k))
\simeq S_k/(I(G_k) +(x_1^2, \ldots, x_k^2))
\]
is Cohen-Macaulay.
Since the polarization preserves Cohen-Macaulayness,
\[
S/(I(G_k) +
(x_1^2, \ldots, x_k^2))^{\mbox{\tiny pol}}
= S/(I(G_k)+ (x_1y_1, \dots, x_ky_k)) = S/I(O_T(G))
\]
is Cohen-Macaulay,
where $(x_1^2, \ldots, x_k^2)^{\mbox{\tiny pol}}$ stands
for the polarization of $(x_1^2, \ldots, x_k^2)$.
See \cite{SV} for basic properties of polarization.

(2) $\Longrightarrow$ (3). Every Cohen-Macaulay ideal is unmixed
\cite{BH}.

(3) $\Longrightarrow$ (1). Suppose $G$ is not Cohen-Macaulay.
We want to prove that there exists a subset $T\subset [n]$ such that
$O_T(G)$ is not unmixed.
Since $G$ is not Cohen-Macaulay the sequence $\{x_i-y_i:
1\leq i\leq n\}$ is not a regular sequence of $S/I(G)$.
Hence there exists $k\geq 1$
such that $\{x_i-y_i: i\in [k-1]\}$ is a regular sequence
of $S/I(G)$ and $x_k-y_k$ is not regular on the ring
\[
R:=S_{k-1}/(I(G_{k-1}) +(x_1^2, \ldots, x_{k-1}^2)) \simeq S/(I(G) + (x_1 -y_1, \ldots, x_{k-1}-y_{k-1})).
\]
Set $J=I(G_{k-1}) +(x_1^2, \ldots, x_{k-1}^2)$.
Since $x_k-y_k$ is not regular on $R$, then
\[
 x_k-y_k\in \bigcup_{P\in \Ass R} P
\]
and there exists an associated prime ideal $\widetilde{P}$ of $J$
such that
$ x_k-y_k \in \widetilde{P}$.
Since $x_k\in \widetilde{P}$ or $y_k\in \widetilde{P}$, we have $x_k, y_k \in \widetilde{P}$.
Hence $\height \widetilde{P}>n$.
Hence $R$ is not unmixed.
Therefore
$S/(I(G_{k-1}) +(x_1^2, \ldots, x_{k-1}^2))^{\mbox{\tiny pol}} \simeq  S/I(O_{T_{k-1}}(G))$ is not unmixed.
%
%
\end{proof}

For distinct $i_{1}, i_{2}, \ldots ,i_{r} \in [n] $
we denote by $C_{i_{1}i_{2} \dots i_{r}}$ the cycle $C$ with
$$
V(C)=\{x_{i_1}, y_{i_1}, x_{i_2}, \ldots , x_{i_r}, y_{i_r}\}
$$
and
$$
E(C)=\{x_{i_1} y_{i_1}, y_{i_1} x_{i_2}, x_{i_2} y_{i_2}, \ldots, y_{i_r} x_{i_r},  y_{i_r} x_{i_1} \}.
$$

\begin{Proposition}\label{Prop:perfect}
Let $G$ be an unmixed graph with $2n$ vertices, which are not isolated,  and
$\height I(G) = n$.
We assume the condition (*).
Then the following conditions are equivalent:
\begin{enumerate}
 \item The subset $\{x_1 y_1, x_2 y_2, \dots, x_n y_n\}$ of $E(G)$
 is a unique perfect matching in $G$.
 \item  The  cycle $C_{ij}$ is not included in G for any $i < j$.
 \item For any $r\ge 2$  the cycle $C_{i_{1}i_{2} \dots i_{r}}$
 is not included in G for any subset $\{i_1, i_2, \dots , i_r \} \subset [n]$ of cardinality $r$.
\end{enumerate}
\end{Proposition}

\begin{proof}

(1) $\Longrightarrow$ (2). Suppose  $C_{ij} \subset G$.
Then we have two perfect matchings in $G$:
\[
\{x_1 y_1,x_2 y_2,\ldots,x_n y_n\},
\]
\[
\{x_1 y_1,x_2 y_2,\ldots,x_{i-1}y_{i-1}, x_iy_j, x_jy_i,
x_{i+1}y_{i+1}, \ldots, x_n y_n\}.
\]

(2)$\Longrightarrow$ (3). We proceed by induction on $r$.

For $r=2$  there is nothing to prove.
Assume $r >2$ and suppose that $C_{i_{1}i_{2} \dots i_{r}} \subset G$.
Since $y_{i_{r-1}}x_{i_r}, y_{i_r}x_{i_{1}} \in E(G)$,
we have $y_{i_{r-1}}x_{i_{1}} \in E(G)$
by Theorem \ref{th:unmix_ijk}.
Hence $C_{i_{1}i_{2} \dots i_{r-1}}  \subset G$,
which is a contradiction with the inductive hypothesis.

(3) $\Longrightarrow$ (1). Suppose there exists another perfect
matching:
\[
\{x_1 y_{i_1},x_2 y_{i_2},\ldots,x_n y_{i_n}\}\subset E(G).
\]
Then we define a permutation $\sigma$ by
\[
\sigma=
\left(
\begin{array}{cccc}
1   & 2   & \ldots & n\\
i_1 & i_2 & \ldots  & i_n
\end{array}
\right).
\]
Then $\sigma$ can be decomposed as  $\sigma=\prod \sigma_i$,
where each $\sigma_i$ is a cycle of $\sigma$.
Since $\sigma$ is not an identity permutation,
for some $i$ the cycle $\sigma_i$ is of the form $(j_1 \,j_2\,\ldots\, j_r)$ with $r\ge 2$.
Then we have
$C_{j_{r}j_{r-1} \dots j_{1}}\subset G$.
\end{proof}

%
%

\begin{Theorem}\label{th:CM}
Let $G$ be an unmixed graph with $2n$ vertices, which are not isolated,  with $\height I(G) = n$ satisfying
the condition (*).
Then the following conditions are equivalent:

\begin{enumerate}
 \item $G$ is Cohen-Macaulay.
 \item $\Delta(G)$ is strongly connected.
 \item  The  cycle $C_{ij}$ is not included in $G$ for any $i < j$.
\end{enumerate}
\end{Theorem}

\begin{proof}
 (1) $\Longrightarrow$ (2). Well known.

 (2) $\Longrightarrow$ (3).
Assume that
$C_{ij} \subset G$
for some $i< j$.
Let $F$ be a facet of $\Delta(G)$ such that $x_i \in F$.
Since $x_i y_j\in E(G)$, we have  $y_j\notin F$ and by
unmixedness of $G$ it follows that $x_j\in F$.
Hence $\{x_i, x_j\} \subset  F$.
Let $F'$ be a facet of $\Delta(G)$ such that $\{y_i, y_j\} \subset F'$.

We show that there does not exist a chain of facets of $\Delta(G)$ such
that
\[
 F=F_0, F_1, \ldots, F_m =F',\, \mbox{ with } \sharp (F_i \cap F_{i+1}) = n-1 \mbox{ for  }i=1, \dots ,m-1.
\]
Every facet $H \in \Delta(G)$ is one of the following form:
\[H = \{z_1, \ldots, z_{i-1}, x_i, z_{i+1},\ldots, z_{j-1},x_j, z_{j+1},\ldots,z_n\}\]
or
\[H = \{ z_1, \ldots, z_{i-1},y_i, z_{i+1}, \ldots, z_{j-1}, y_j, z_{j+1}, \ldots, z_n\},\]
where $z_{k} \in \{ x_k, y_k\}$,
since $\{x_i y_i, x_j y_j, x_i y_j, x_j y_i\} \subset E(G)$.
Hence it is impossible to find such a chain.
Hence $\Delta(G)$ is not strongly connected.

(3) $\Longrightarrow$ (1). In order to prove  the statement  by Proposition \ref{Prop:CM}
it is sufficient to verify that $O_T(G)$ is unmixed for every subset $T$
of $[n]$.
Hence we prove that conditions (i) and (ii)  of Proposition \ref{th:unmix_ijk} are satisfied
for the graph $G'=O_T(G)$.

First we check the condition (i)  for $G'$.
We may assume that $j \notin T$.

Suppose $i \notin T$. We must show the following:
\newline
``If $z_i x_j$, $y_j x_k\in E(G')$, then $z_i x_k\in E(G')$ for
distict $i$, $j$ and $k$ and for $z_i\in \{ x_i, y_i\}$."
\newline
Since $z_i x_j$, $y_j x_k\in E(G)$ and $G$ is unmixed, by Theorem \ref{th:unmix_ijk}
we have $z_i x_k\in E(G)$.
Hence $z_i x_k\in E(G')$.

Suppose $i \in T$. We must show the following:
\newline
``If $x_i x_j$, $y_j x_k\in E(G')$ then $x_i x_k\in E(G')$
for distict $i$, $j$ and $k$."
\newline
Since $y_j x_k\in E(G')$, we have $y_j x_k\in E(G)$. Since $x_i
x_j \in E(G')$, we have either $x_i x_j\in E(G)$ or $y_i x_j\in
E(G)$. If $x_i x_j\in E(G)$, then by Theorem \ref{th:unmix_ijk}
we have $x_i x_k\in E(G)$, since $y_j x_k\in E(G)$ and $G$ is
unmixed. Similarly,  if $y_i x_j\in E(G)$, then we have $y_i
x_k\in E(G)$, since $y_j x_k\in E(G)$. In both cases, we have $x_i
x_k\in E(G')$.

%
%
%
%
%
Next we check the condition (ii) for $G'$. We may assume that  $j
\notin T$. We also assume that  $i \in T$. We must show that
either $x_i x_j \notin E(G')$ or $x_i y_j \notin E(G')$. Suppose
$x_i x_j , x_i y_j \in E(G')$. Then we have $x_i y_j \in E(G)$,
and either $x_i x_j \in E(G)$ or $y_i x_j \in E(G)$. Since $G$ is
unmixed, $x_i x_j \in E(G)$ is impossible by  Theorem
\ref{th:unmix_ijk}, (ii). While the condition  $y_i x_j \in E(G)$
is also impossible, since $G$ does not have the cycle $C_{ij}$
for any $i < j$. It is a contradiction.
\end{proof}

The next lemma is crucial for giving another criterion for the
Cohen-Macaulayness of our graphs.

\begin{Lemma}\label{lem:order} Let $G$ be an unmixed graph
with $2n$ vertices, which are not isolated,  and $\height I(G) =n$. We assume the
condition (*).

If $G$ is a Cohen-Macaulay graph then there exists a suitable simultaneous change
of labeling on both $\{x_i\}$ and $\{y_i\}$
(i.e.,
we relable  $(x_{i_1}, \ldots, x_{i_n})$ and $(y_{i_1}, \ldots,
y_{i_n})$ as  $(x_1, \ldots, x_n)$ and $(y_1, \ldots, y_n)$ at
the same time), such that $x_iy_j \in E(G)$ implies $i \leq j$.
\end{Lemma}

\begin{proof}
We can define a partial order $\preceq$ on $X$
by
\[
x_i \preceq x_j  \mbox{ if and only if } x_i y_j \in E(G).
\]
In fact, the reflexivity  holds by (*), the transitivity holds by
unmixedness of $G$ (see Theorem 2.4 (i)) and the antisymmetry holds
since $G$ contains no cycle $C_{ij}$ for any $i < j$. Take a linear
extension of $\preceq$, which we call $\preceq '$. By the linear
order $\preceq '$, we have $x_{i_1} \preceq '  \cdots \preceq '
x_{i_n}$. We relabel them as $x_{1} \preceq '  \cdots \preceq '
x_{n}$. At the same time we relabel $y_{i_1},\ldots , y_{i_n}$ as
$y_{1},\ldots , y_{n}$. Then if $x_iy_j \in E(G)$, $x_{i} \preceq '
x_{j}$. Hence $i \leq j$.
%
%
%
%
%
%
%
\end{proof}

Hence for a Cohen-Macaulay graph $G$ with $2n$ vertices, which are not isolated,
and $\height I(G) =n$ satisfying the condition (*), we may assume
that

\vspace{.2in} (**) $x_iy_j \in E(G)$ implies $i \leq j$.
\vspace{.2in}

Now we state another Cohen-Macalay criterion on our graphs, which is
generalization of Herzog-Hibi (\cite {HH}, Theorem 3.4).

\begin{Theorem}\label{Crit:CM} Let $G$ be a graph with
$2n$ vertices, which are not isolated,  and with $\height I(G) =n$.
We assume the conditions (*) and (**).

Then the following conditions are equivalent:
\begin{enumerate}
\item $G$ is Cohen-Macaulay;
\item $G$ is unmixed;
\item The following conditions hold:
\begin{enumerate}
\item[(i)] If $z_ix_j, y_jx_k \in E(G)$ then $z_ix_k \in E(G)$ for distict $i, j, k$ and
for $z_i \in \{ x_i, y_i \}$;
\item[(ii)] If $x_iy_j \in E(G)$ then $x_ix_j \notin E(G)$.
\end{enumerate}
\end{enumerate}
\end{Theorem}
\begin{proof}
(1) $\Longrightarrow$ (2) is well known.\newline
(2) $\Longrightarrow$ (1) follows from Theorem \ref{th:CM}, since we assume the condition (**).\newline
(2) $\Longleftrightarrow$ (3) follows from Theorem \ref{th:unmix_ijk}.
\end{proof}

As an easy consequence of the previous results we obtain the upper bound for the minimal number $\mu (I(G))$
of generators of $I(G)$:

\begin{Corollary}
Let $G$ be a graph with
$2n$ vertices, which are not isolated,  and with $\height I(G) =n$.
\begin{enumerate}
\item[(i)] If $G$ is unmixed, then $\mu(I(G)) \leq n^2$.
\item[(ii)] If $G$ is Cohen-Macaulay, then $\mu(I(G)) \leq \frac{n(n+1)}2$.
\end{enumerate}
\end{Corollary}

\begin{proof} The statements are consequences of the criteria for
the unmixedness and for the Cohen-Macaulayness given by Proposition
\ref{th:unmix_ijk} and Theorem \ref{Crit:CM}.
\end{proof}

\section{Shellability and Cohen-Macaulay type}

In this section if $G$ is a graph such that $\sharp V(G) = 2n$ and
$\height I(G) = n$, we show the equivalence between Cohen-Macaulayness of $G$ and
shellability of the complementary simplicial complex $\Delta(G)$.
We also express the Cohen-Macaulay type of $S/I(G)$ in a combinatorial way.

\begin{Theorem} \label{th:shell} Let $G$ be an unmixed graph
with $2n$ vertices, which are not isolated,  and with $\height I(G) =n$.
Then $G$ is Cohen-Macaulay if and only if $\Delta(G)$ is shellable.
\end{Theorem}

We just give a proof of the following lemma. The rest of the proof
is almost identical with the proof of (\cite{EV}, Theorem 2.9).

\begin{Lemma} \label{lem:deg}
Let $G$ be a Cohen-Macaulay graph with $2n$ vertices, which are not isolated,
and $\height I(G)=n$.
Then there exists a vertex $v\in V(G)$ such that $\deg(v)=1$.
\end{Lemma}

\begin{proof}
We may assume the condition (*).
Suppose that each $v\in V(G)$ has at least degree 2.
Let $i_1, i_2, \ldots $ be  a  sequence  such that
$y_{i_1} x_{i_2}, y_{i_2} x_{i_3}, \ldots \in E(G)$
with $i_j\neq i_{j+1}$.
Since the cardinality of $Y$ is finite,
there must be exist integers $s<t$
such that $i_t=i_s$.
We may assume that $i_s,i_{s+1},\ldots,i_{t-1}$  are distinct.
This induces the cycle $C_{i_{s}i_{s+1} \cdots i_{t-1}}\subset G$.
Therefore $G$ is not Cohen-Macaulay by Proposition \ref{Prop:perfect} and Theorem \ref{th:CM}.
\end{proof}

Now we express the Cohen-Macaulay type of a graph belonging to
our class, imitating the bipartite case (see \cite {Vi1}, pp.
184-185).

\begin{Lemma} \label{lem:soc} Let $G$ be a Cohen-Macaulay graph with
$2n$ vertices, which are not isolated,  and $\height I(G) =n$.
We assume the condition (*).
Then
\[
\Soc \ (K[x_1, \ldots, x_n]/\left(I(O_{[n]}(G)\vert_{X}) + (x_1^2, \ldots, x_n^2)\right))
\]
is generated by all the monomials $x_{i_1} \cdots x_{i_r}$
such that $\{x_{i_1}, \ldots , x_{i_r}\}$ is a maximal independent set of $O_{[n]}(G)\vert_{X}$.
\end{Lemma}

\begin{proof}
The ring $A:=K[x_1, \ldots, x_n]/\left(I(O_{[n]}(G)\vert_{ X}) + (x_1^2, \ldots, x_n^2)\right)$
is spanned as  a $K$-vector space  by the image of $1$
and the images of the squarefree monomials
\begin{equation}\label{monomials}
    x_{i_1}\cdots  x_{i_r},\qquad 1 \leq i_1 < i_2 < \ldots < i_r
    \leq n
\end{equation}
such that $x_{i_j}x_{i_k} \notin E(O_{[n]}(G)\vert_{X})$, for $j
\neq k$, i.e. $\{x_{i_1},\ldots, x_{i_r}\}$ is an independent set of
$O_{[n]}(G)\vert_{X}$.
Since $A$ is an artinian positively graded algebra, Soc $A = (0:_A A_+)$ is generated
by the images of the squarefree monomials of the form (\ref{monomials}) such that
$\{x_{i_1},\ldots, x_{i_r}\}$ is a maximal independent set of
$O_{[n]}(G)\vert_{ X}$.
\end{proof}

\begin{Corollary}
Let $G$ be a Cohen-Macaulay graph
with $2n$ vertices, which are not isolated,  and $\height I(G) =n$.
We assume the condition (*).
Then
\begin{enumerate}
\item[(i)]
$\mbox{\rm type } S/I(G) = \sharp \Upsilon(O_{[n]}(G)\vert_{X})$,
where $\Upsilon(O_{[n]}(G)\vert_{X})$ is the family of all minimal vertex covers of $O_{[n]}(G)\vert_{X}$.
In particular, {\rm type }$ S/I(G)$ is independent from the base field $K$.
\item[(ii)]
$G$ is level if and only if $O_{[n]}(G)\vert_{X}$ is unmixed.
In particular, level-ness of $G$ is independent from the base field $K$.
\end{enumerate}
\end{Corollary}

\begin{proof}
Set $S = K[x_1, \ldots, x_n,y_1, \ldots, y_n]$ and $S_n = K[x_1, \ldots, x_n]$.

(i)Since $G$ is Cohen-Macaulay and $\{x_1 - y_1, \ldots, x_n -y_n\}$ is a regular sequence,
we have
\begin{eqnarray*}
\mbox{\rm type } S/I(G)
& = & \dim _K \mbox{\rm Soc } S/\left(I(G) + (x_1 - y_1, \ldots, x_n -y_n)\right) \\
& = & \dim _K \mbox{\rm Soc } S_n/\left(I(O_{[n]}(G)\vert_{X}) +(x_1^2, \ldots, x_n^2)\right) \\
& = & \sharp \Upsilon (O_{[n]}(G)\vert_{X})
\end{eqnarray*}
by the previous lemma.

(ii) When $G$ is Cohen-Macaulay, $G$ is level if and only if
$$
\Soc  S/\left(I(G) + (x_1 - y_1, \ldots, x_n -y_n)\right)
$$ is equi-generated.
By the previous lemma it is equivalent to that $O_{[n]}(G)\vert_{X}$ is unmixed.
\end{proof}

\begin{Corollary}
Let $G$ be a Cohen-Macaulay graph
with $2n$ vertices, which are not isolated,  and $\height I(G) =n$.
We assume the condition (*).
Then the following conditions are equivalent:
\begin{enumerate}
\item $G$ is Gorestein;
\item $I(G) = (x_1y_1, \ldots, x_ny_n)$;
\item $G$ is a complete intersection.
\end{enumerate}
\end{Corollary}

\begin{proof} (1) $\Rightarrow$ (2). $G$ is Gorenstein if and only
if $S/I(G)$ is Cohen-Macaulay and $\mbox{type } S/I(G) = 1$.
Since $1= \mbox{type } S/I(G) = \sharp \Upsilon(O_{[n]}(G)\vert_{X})$, it
follows that $O_{[n]}(G)\vert_{X}$ has a unique
minimal vertex cover.
Hence $O_{[n]}(G)\vert_{X}$ is isolated $n$ vertices.
Hence $I(G) = (x_1y_1, \ldots, x_ny_n)$.

(2) $\Rightarrow$ (3). From its definition.

(3) $\Rightarrow$ (1). See \cite{BH}.
\end{proof}

\section{B-grafted graph}

In this section we introduce a new class of graphs $G$ with $\sharp V(G) = 2n$ and
with $\height I(G) = n$ and we study its Cohen-Macaulayness.

Let $H_0$ be a graph with the labeled vertices $ 1,2, \ldots , p$.

For every $i = 1, \ldots, p$ let $B_i$ be a bipartite graph with
labeled partition $X_i$ and $Y_i$ such that $\sharp X_i = \sharp
Y_i = n_i$. (We do not give a label to each vertex of $B_i$, but
we distinguish the partition $X_i$ and $Y_i$.) We assume that
$B_i$ has no isolated vertex for every $i = 1, \ldots, p$. We
define the graph
\[G = G(H_0; B_1, \ldots, B_p)\] as follows:
The vertex set of $G$ is $V(G) := X\cup Y$, where $X=X_1 \cup \ldots \cup X_p$, and $Y= Y_1 \cup
\ldots \cup Y_p$.
The edge set $E(G)$ of $G$ defined by:
\[xy \in E(G)\mbox{ if and only if }\]
either
\[\mbox{there exist}\,\,\, i, j\,\,\,\mbox{ such that }\,\,\,\,x\in X_i, y \in X_j, \mbox{ and  }ij \in E(H_0)\]
or
\[\mbox{there exists}\,\,\, i\,\,\,
\mbox{ such that }\,\,\,\,x\in X_i, y \in Y_i, \mbox{ and } xy \in E(B_i).\]
We call such a graph $G$ the $B$-\textit{grafted graph}.
Note that $X$ is a minimal vertex cover of $G$ and $Y$ is a maximal independent set of $G$.
Note also that $\sharp V(G) = 2(\sum_{i=1}^p n_i)$.

\medskip

\begin{Example}\label{ex:bgraft}\rm
Let $H_0$ be a cycle of the length 3. By the following bipartite
graphs $B_1$, $B_2$, $B_3$,  we obtain the $B$-grafted graph $G$:
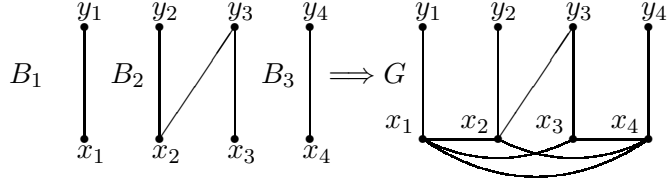
\[
\setlength{\unitlength}{.5cm}
\begin{picture}(15,5)

\put(1,1){\line(0,1){3}}

\put(3,1){\line(0,1){3}}
\put(3,1){\line(2,3){2}}

\put(5,1){\line(0,1){3}}

\put(7,1){\line(0,1){3}}

\put(1,1){\circle*{.2}}
\put(.8,.5){$x_1$}

\put(3,1){\circle*{.2}}
\put(2.8,.5){$x_2$}

\put(5,1){\circle*{.2}}
\put(4.8,.5){$x_3$}

\put(7,1){\circle*{.2}}
\put(6.8,.5){$x_4$}

\put(1,4){\circle*{.2}}
\put(.8,4.3){$y_1$}

\put(3,4){\circle*{.2}}
\put(2.8,4.3){$y_2$}

\put(5,4){\circle*{.2}}
\put(4.8,4.3){$y_3$}

\put(7,4){\circle*{.2}} \put(6.8,4.3){$y_4$}

\put(-1,2.5){$B_1$}

\put(1.7,2.5){$B_2$}

\put(5.7,2.5){$B_3$}

\put(7.5,2.5){$\Longrightarrow G$}

\qbezier(10,1)(12,0)(14,1)

\qbezier(12,1)(14,0)(16,1)

\qbezier(10,1)(13,-1)(16,1)

\put(10,1){\line(0,1){3}}
\put(10,1){\line(1,0){2}}

\put(12,1){\line(0,1){3}}
\put(12,1){\line(2,3){2}}

\put(14,1){\line(0,1){3}}
\put(14,1){\line(1,0){2}}

\put(16,1){\line(0,1){3}}

\put(10,1){\circle*{.2}}
\put(9,1.3){$x_1$}

\put(12,1){\circle*{.2}}
\put(11,1.3){$x_2$}

\put(14,1){\circle*{.2}}
\put(13,1.3){$x_3$}

\put(16,1){\circle*{.2}}
\put(15,1.3){$x_4$}

\put(10,4){\circle*{.2}}
\put(9.8,4.3){$y_1$}

\put(12,4){\circle*{.2}}
\put(11.8,4.3){$y_2$}

\put(14,4){\circle*{.2}}
\put(13.8,4.3){$y_3$}

\put(16,4){\circle*{.2}}
\put(15.8,4.3){$y_4$}


\end{picture}
\]
\end{Example}

\medskip

\begin{Remark} \rm If $B_i$ is just a complete graph with 2 vertices, i.e.,
a complete bipartite graph with $\sharp X_i = \sharp Y_i = 1$ for $i= 1, \ldots, p$,
then the B-grafted graph $G$ is called a grafted graph in \cite{SF}.
\end{Remark}

\begin{Theorem} \label{B-graph:ht}
The $B$-grafted graph $G(H_0; B_1, \ldots, B_p)$ is
Cohen-Macaulay (unmixed, respectively) if and only if every
bipartite graph $B_i$ is Cohen-Macaulay (unmixed, respectively)
for $i = 1, \ldots,p$.

\end{Theorem}
\begin{proof}
It is clear from Theorem \ref{th:CM} (Proposition
\ref{th:unmix_ijk}, respectively).
\end{proof}
\medskip\noindent
{\bf Acknowledgments.}
The third author acknowledges the financial support of GNSAGA-INDAM and the hospitality
during his stay at the  Department of Mathematics of the University of Messina (Italy).
This work was also supported by KAKENHI18540041 and KAKENHI20540047.

\bibliographystyle{plain}

\end{document}